\batchmode
\makeatletter
\def\input@path{{\string"/Users/paranoia/Documents/Research/mypapers/An EL-labeling of L(G)/\string"/}}
\makeatother
\documentclass[12pt,oneside]{amsart}
\usepackage[T1]{fontenc}
\usepackage[latin9]{inputenc}
\usepackage{verbatim}
\usepackage{amssymb}

\makeatletter
\numberwithin{equation}{section} 
\numberwithin{figure}{section} 
  \theoremstyle{plain}
  \newtheorem{thm}{Theorem}[section]
  \theoremstyle{remark}
  \newtheorem{note}[thm]{Note}
  \theoremstyle{plain}
  \newtheorem{conjecture}[thm]{Conjecture}
 \theoremstyle{definition}
  \newtheorem{example}[thm]{Example}
  \theoremstyle{plain}
  \newtheorem{lem}[thm]{Lemma}
  \theoremstyle{plain}
  \newtheorem{cor}[thm]{Corollary}
  \theoremstyle{plain}
  \newtheorem{prop}[thm]{Proposition}

\makeatother

\begin{document}
%
{}

\newcommand{\normalin}{\mathrel{\triangleleft}}

\newcommand{\innormal}{\mathrel{\triangleright}}

\newcommand{\semidirect}{\mathbin{\rtimes}}

\newcommand{\Stab}{\operatorname{Stab}}

%
{}

\newcommand{\bdry}{\partial}

\newcommand{\susp}{\operatorname{susp}}

%
{}

\newcommand{\lrprod}{\mathop{\check{\prod}}}

\newcommand{\lrtimes}{\mathbin{\check{\times}}}

\newcommand{\urtimes}{\mathbin{\hat{\times}}}

\newcommand{\urprod}{\mathop{\hat{\prod}}}

\newcommand{\subsetdot}{\mathrel{\subset\!\!\!\!{\cdot}\,}}

\newcommand{\dotsupset}{\mathrel{\supset\!\!\!\!\!\cdot\,\,}}

\newcommand{\precdot}{\mathrel{\prec\!\!\!\cdot\,}}

\newcommand{\dotsucc}{\mathrel{\cdot\!\!\!\succ}}

\newcommand{\des}{\operatorname{des}}

%
{}

\newcommand{\modreln}{\mathrel{M}}

%
{}

\newcommand{\link}{\operatorname{link}}

\newcommand{\freejoin}{\mathbin{\circledast}}

\newcommand{\stellarsd}{\operatorname{stellar}}

\newcommand{\conv}{\operatorname{conv}}

%
{}

\newcommand{\cosetposet}{\overline{\mathfrak{C}}}

\newcommand{\cosetlat}{\mathfrak{C}}

\newcommand{\nlat}{\mathcal{N}}

\newcommand{\nproj}{\rho}


\title{An $EL$-labeling of the subgroup lattice}

\author{Russ Woodroofe}

\email{russw@math.wustl.edu}

\begin{abstract}
In a 2001 paper, Shareshian conjectured that the subgroup lattice
of a finite, solvable group has an $EL$-labeling. We construct such
a labeling and verify that our labeling has the expected properties.
\end{abstract}
\maketitle

\section{Introduction}

All groups, posets, and simplicial complexes in this paper are finite.
We recall that the \emph{subgroup lattice} $L(G)$ of a group $G$
is the set of all subgroups of the group, ordered under inclusion.
$L(G)$ is a lattice, with $H\wedge K=H\cap K$ and $H\vee K=\langle H,K\rangle$. 

Any poset $P$ is closely associated with its \emph{order complex}
$\vert P\vert$, a simplicial complex with faces the chains in $P$.
Considering the order complex allows us to use combinatorial topology
definitions and theorems with $P$. One such definition is that of
{}``shellability.'' A {}``shellable'' complex is essentially one
where the facets fit nicely together \cite{Bjorner:1980,Bjorner/Wachs:1983,Bjorner/Wachs:1996,Bjorner/Wachs:1997};
the precise definition will not be important to us. A shellable poset
is one with shellable order complex.

The connection with the subgroup lattice is surprising and beautiful:

\begin{thm}
\emph{\label{thm:Shareshian-IffVersion}(Shareshian \cite[Theorem 1.4]{Shareshian:2001})}
$L(G)$ is shellable if and only if $G$ is solvable.
\end{thm}
Let us talk about the techniques used to prove the {}``if'' direction
of Theorem \ref{thm:Shareshian-IffVersion}. There are two main techniques
to show that a bounded poset is shellable, both developed by Björner
and Wachs \cite{Bjorner:1980,Bjorner/Wachs:1983,Bjorner/Wachs:1996,Bjorner/Wachs:1997}.
The first is to label the edges of the Hasse diagram in a manner such
that on every interval:

\begin{enumerate}
\item There is a unique chain where the labels (read from bottom to top)
are increasing.
\item The unique increasing chain is lexicographically first.
\end{enumerate}
A labeling satisfying these two properties is called an \emph{$EL$-labeling}.

The second is to label the atoms of the poset. A \emph{recursive atom
ordering} of a bounded poset $P$ is an ordering $a_{1},a_{2},\dots$
of the atoms of $P$ such that

\begin{enumerate}
\item For any $j$, the interval $[a_{j},\hat{1}]$ has a recursive atom
ordering in which the atoms in $[a_{j},\hat{1}]$ that are above some
$a_{i}$ for $i<j$ come first.
\item For all $i<j$, and $x$ with $a_{i},a_{j}<x$, there is a $k<j$
and an atom $y<x$ of $[a_{j},\hat{1}]$ with $a_{k}<y$.
\end{enumerate}
A bounded poset with either an $EL$-labeling or a recursive atom
ordering is shellable. The two are somewhat related: a poset with
a recursive atom ordering has a {}``$CL$-labeling'', which is a
generalization of the idea of an $EL$-labeling. As a poset is shellable
if and only if its dual is shellable, recursive coatom orderings and
dual $EL$-labelings are also of interest.

Shareshian proved the {}``if'' direction of Theorem \ref{thm:Shareshian-IffVersion}
as follows:

\begin{thm}
\emph{(Shareshian \cite[Corollary 4.10]{Shareshian:2001})} If $G$
is solvable, then $L(G)$ has a recursive coatom ordering.
\end{thm}
\begin{note}
Interestingly, the ordering of maximal subgroups (coatoms) that Shareshian
used had already been studied by Doerk and Hawkes \cite[Chapter A.16, especially Definition 16.5]{Doerk/Hawkes:1992}.
\end{note}
An $EL$-labeling gives useful information about a poset. For example,
one of the nicest consequences is that the set of descending chains
forms a cohomology basis for $\vert P\vert$. Unfortunately, although
every poset with a recursive (co-)atom ordering has a (dual) $CL$-labeling,
the construction is complicated enough that nice enumerative results
(such as the cohomology basis) coming from $EL$/$CL$-labelings are
usually difficult or impossible to use.

\medskip{}
The topology of the subgroup lattice of a solvable group had been
studied before Shareshian. Let $G$ be a solvable group, with chief
series $1=N_{0}\subset N_{1}\subset\dots\subset N_{k}=G$. A \emph{complement}
to a subgroup $N$ is a subgroup $H$ with $HN=G$ and $H\cap N=1$.
A \emph{chain of complements} (to the given chief series) is a chain
$1=H_{k}\subset H_{k-1}\subset\dots\subset H_{0}=G$ where $H_{i}$
is a complement to $N_{i}$ (for each $i$). Then

\begin{thm}
\emph{\label{thm:SubgroupLatticeTopology}(Thévenaz \cite[Theorem 1.4]{Thevenaz:1985})}
For any solvable group $G$, $\vert L(G)\vert$ has the homotopy type
of a wedge of spheres of dimension $k-2$, and the spheres are in
bijective correspondence with the chains of complements to any given
chief series.
\end{thm}
In light of the cohomology basis mentioned above, Theorem \ref{thm:SubgroupLatticeTopology}
naturally leads to the following conjecture:

\begin{conjecture}
\emph{\label{con:SglatticeHasELlabeling}(Shareshian \cite[Conjecture 1.6]{Shareshian:2001})}
For any solvable group $G$, $L(G)$ admits an $EL$-labeling where
the descending chains are the chains of complements to a chief series.
\end{conjecture}
In the rest of this paper, we will extend the theory of left modular
lattices to construct both an $EL$-labeling and a dual $EL$-labeling
satisfying Conjecture \ref{con:SglatticeHasELlabeling}.

\section{Left modularity}

Our starting point will be left modularity. Let $L$ be any lattice.
An element $x\in L$ is \emph{left modular} if for all $y<z$ we have
$(y\vee x)\wedge z=y\vee(x\wedge z)$, i.e., if it satisfies one side
of the requirement for modularity. 

\begin{example}
The Dedekind identity (see for example \cite[1.3.14]{Robinson:1996})
says that $H(N\cap K)=HN\cap K$ for any subgroup $N$, and subgroups
$H\subseteq K$ of a group $G$. Since a normal subgroup $N$ of $G$
satisfies $HN=NH=\langle H,N\rangle$ for every subgroup $H$, a normal
subgroup is left modular in $L(G)$.

\medskip{}
Liu gave a helpful alternative characterization of left modular elements.
Let $y\lessdot z$ denote a \emph{cover relation}, that is, a pair
$y<z$ such that if there is an $x$ with $y\leq x\leq z$, then $x=y$
or $x=z$.
\end{example}
\begin{thm}
\emph{\label{thm:LeftModularAltChar}(Liu} \cite[Theorem 2.1.4]{Liu:1999}\emph{,
also in \cite[Theorem 1.4]{Liu/Sagan:2000})} Let $x$ be an element
in a lattice $L$. The following are equivalent.
\begin{enumerate}
\item $x$ is left modular.
\item For any $y<z$ we have $x\vee z\neq x\vee y$ or $x\wedge z\neq x\wedge y$.
\item For any $y\lessdot z$ we have $x\vee z=x\vee y$ or $x\wedge z=x\wedge y$,
but not both.
\end{enumerate}
\end{thm}
Part (3) of Theorem \ref{thm:LeftModularAltChar} leads us to the
following definition: let $\hat{0}=x_{0}<x_{1}<\dots<x_{k}=\hat{1}$
be a (not necessarily maximal) chain with every $x_{i}$ left modular.
Then we say $x_{i+1}/x_{i}$ \emph{weakly separates} a cover relation
$y\lessdot z$ if $x_{i}\wedge z=x_{i}\wedge y$ but $x_{i+1}\vee z=x_{i+1}\vee y$.
Any given cover relation is weakly separated by a unique $x_{i+1}/x_{i}$
in the modular chain.

Then it is natural to consider the labeling \[
\lambda(y\lessdot z)=i\quad\quad\mbox{where }x_{i+1}/x_{i}\mbox{ weakly separates }y\lessdot z\mbox{.}\]

\begin{thm}
\emph{(Liu \cite[Theorem 3.2.6]{Liu:1999})} If the left modular chain
$\hat{0}=x_{0}<x_{1}<\dots<x_{k}=\hat{1}$ is a maximal chain, then
$\lambda$ is an $EL$-labeling.
\end{thm}
In this situation (where $L$ has a maximal chain of left modular
elements) we say that $L$ is \emph{left modular}. Left modular lattices
have been studied in several papers \cite{McNamara/Thomas:2006,Thomas:2005,Blass/Sagan:1997}
in addition to the ones already referenced. Lattices with chains of
modular elements were studied in \cite{Hersh/Shareshian:2006}.

\medskip{}
Motivated by the situation in a solvable group (where the chief series
is a left modular chain, but not necessarily a maximal one), we ask
what happens with the labeling $\lambda$ when $\hat{0}=x_{0}<x_{1}<\dots<x_{k}=\hat{1}$
is \underbar{not} maximal. We don't get an $EL$-labeling, but we
can still say some things about the increasing chains on an interval. 

Let $[w,z]$ be an interval in $L$. Then $w\leq w\vee x_{i}\wedge z\leq z$
for all $i$, and we notice that $w\lneq w\vee x_{i}\wedge z$ for
large enough $i$ (in particular, $i=k$ gives $w\vee\hat{1}\wedge z=z$).
So let $c_{0}=w$, and inductively construct $c_{j}$ as follows:
let $i(j)$ be the maximal index such that $c_{j}\vee x_{i(j)}\wedge z=c_{j}$.
Then let \[
c_{j+1}=c_{j}\vee x_{i(j)+1}\wedge z=w\vee x_{i(j)+1}\wedge z.\]
This gives a chain $\mathbf{c}=\{w=c_{0}<c_{1}<\dots<c_{m}=z\}$ between
$w$ and $z$. Every edge on the interval $[c_{j},c_{j+1}]$ receives
an $i(j)$ label, since for every $y$ on $[c_{j},c_{j+1}]$ we have
\[
y\vee x_{i(j)+1}=y\vee(x_{i(j)+1}\wedge z)\vee x_{i(j)+1}=c_{j+1}\vee x_{i(j)+1};\]
while $y\vee x_{i(j)}\wedge z=y$, so that each $y\vee x_{i(j)}$
is distinct.

\begin{lem}
\label{lem:LMchainGivesPreEL}A maximal chain on $[w,z]$ is (weakly)
increasing if and only if it is an extension of $\mathbf{c}$.
\end{lem}
\begin{proof}
Every extension of $[c_{j},c_{j+1}]$ has every edge labeled with
$i(j)$. Since, by the construction, $i(0)<i(1)<\dots<i(m-1)$, every
maximal extension of $\mathbf{c}$ is (weakly) increasing.

In the other direction, notice that since $w\vee x_{i(0)}\wedge z=w$,
but $w\vee x_{i(0)+1}\wedge z\gneq w$, there must be an edge $d_{j}\lessdot d_{j+1}$
in any maximal chain $\mathbf{d}=\{w=d_{0}<d_{1}<\dots<z\}$ such
that $d_{j}\ngeq w\vee x_{i(0)+1}\wedge z$ but $d_{j+1}\geq w\vee x_{i(0)+1}\wedge z$.
Clearly such an edge receives an $i(0)$ label, and since by the definition
of the labeling any maximal chain cannot have labels less than $i(0)$,
any weakly increasing maximal chain must start with $i(0)$ labels. 

The first edge of $\mathbf{d}$ receives the label $i(0)$ only if
$d_{0}\vee x_{i(0)+1}=d_{1}\vee x_{i(0)+1}$; thus, \[
d_{1}\leq d_{1}\vee x_{i(0)+1}\wedge z=d_{0}\vee x_{i(0)+1}\wedge z=c_{1},\]
and so the first edge of $\mathbf{d}$ is in $[c_{0},c_{1}]$. Repeating
this argument inductively on $[d_{1},z]$ gives that $\mathbf{d}$
is an extension of $\mathbf{c}$, as desired.
\end{proof}
\begin{cor}
\label{cor:LMchainGivesPreEL-lex}A maximal chain on $[w,z]$ is (tied
for) lexicographically first if and only if it is an extension of
$\mathbf{c}$.
\end{cor}
\begin{note}
There is not in general a unique lexicographically first or increasing
chain, as $\mathbf{c}$ may have many extensions.
\end{note}

\begin{note}
We use the term {}``weakly separated'' to highlight that a maximal
chain might have multiple $i$ labels. One might say that $y\lessdot z$
was \emph{separated} by $x_{i+1}/x_{i}$ if the edge was weakly separated
and also $x_{i+1}\wedge y=x_{i}\wedge y$ and $x_{i+1}\vee z=x_{i}\vee z$
(but we will not use this). 
\end{note}
In Section 3, we will show that intervals in $L(G)$ with repeated
$i$ labels are isomorphic to certain sublattices of $[N_{i},N_{i+1}]$,
and in Section 4 we will use this isomorphism to refine $\lambda$
to an $EL$-labeling in the subgroup lattice (of a solvable group).

\section{Projecting into $[N_{i},N_{i+1}]$}

Let $G$ be a solvable group with a chief series $1=N_{0}\subset N_{1}\subset\dots\subset N_{k}=G$,
and let $H$ be any subgroup. The subgroups of $L(G)$ that are normalized
by $H$ form a sublattice $L_{H}(G)$. In this section we will relate
certain sections \[
\nlat_{i}(H)\triangleq[N_{i},N_{i+1}]\cap L_{H}(G)\]
of this lattice to weak separation by the chief series. First:

\begin{lem}
For any $H$, $\nlat_{i}(H)$ is a modular lattice.
\end{lem}
\begin{proof}
$\nlat_{i}(H)$ is closed under intersection and join, so it is a
sublattice of $[N_{i},N_{i+1}]$. By the Correspondence Theorem \cite[1.4.6]{Robinson:1996},
we have that $[N_{i},N_{i+1}]\cong L(N_{i+1}/N_{i})$. Since $N_{i+1}/N_{i}$
is abelian, $[N_{i},N_{i+1}]$ is a modular lattice, and sublattices
of a modular lattice are modular.
\end{proof}
Second, we have a relationship between weak separation of an edge
in $L(G)$ and $\nlat_{i}$.

\begin{lem}
\label{lem:NlatPreserved}If $E\subsetdot F$ is weakly separated
by $N_{i+1}/N_{i}$, then $\nlat_{i}(E)=\nlat_{i}(F)$.
\end{lem}
\begin{proof}
$N_{i+1}E=N_{i+1}F$, so $F\subseteq EN_{i+1}$. Since every subgroup
$N$ in the interval $[N_{i},N_{i+1}]$ is normalized by $N_{i+1}$,
we see that if $E$ normalizes $N$, then so does $F$. The converse
is immediate.
\end{proof}
\begin{note}
When we are looking at an edge or chain(s) of edges that are weakly
separated by $N_{i+1}/N_{i}$, we will often simply write $\nlat_{i}$
to mean $\nlat_{i}(E)=\nlat_{i}(F)=\dots$. Lemma \ref{lem:NlatPreserved}
tells us that this notation makes sense.
\end{note}
Finally, we construct a projection map from $L(G)$ to $[N_{i},N_{i+1}]$.
Let \[
\nproj_{i}(H)=N_{i}\vee H\wedge N_{i+1}=N_{i}H\cap N_{i+1}.\]
It is clear that this is really in $[N_{i},N_{i+1}]$. In fact, $\nproj_{i}(H)$
is in $\nlat_{i}(H)$ (since $N_{i}$, $N_{i+1}$, and $H$ are all
normalized by $H$). Much more is true. Let $[W,Z]_{\mathcal{S}}$
denote the interval $[W,Z]$ in the sublattice $\mathcal{S}$ of $L(G)$;
that is, let $[W,Z]_{\mathcal{S}}$ consist of all $H\in\mathcal{S}$
that are between $W$ and $Z$. 

\begin{prop}
\label{pro:ProjIsomorphism}If there is a chain on the interval $[W,Z]$
with every edge weakly separated by $N_{i+1}/N_{i}$, then $\rho_{i}$
on $[W,Z]$ gives a poset isomorphism \[
[W,Z]_{L(G)}\cong[\nproj_{i}(W),\nproj_{i}(Z)]_{\nlat_{i}}.\]

\end{prop}
\begin{example}
Consider the alternating group on 4 elements with the normal series
$N_{0}=1$, $N_{1}$ the Klein 4 subgroup, and $N_{2}=A_{4}$. Then
$\langle(1\,2\,3)\rangle\subsetdot A_{4}$ is weakly separated by
$N_{1}/N_{0}$, and it projects to $N_{0}\subsetdot_{\nlat_{i}}N_{1}$,
an edge in the sublattice $\nlat_{i}=\nlat_{i}(A_{4})$. Notice that,
although $N_{0}\subset N_{1}$ is a cover relation in $\nlat_{i}$,
it is \underbar{not} a cover relation in $L(G)$, as $N_{0}=1\subset\langle(1\,2)(3\,4)\rangle\subset N_{1}$.
\end{example}
\begin{proof}
(of Proposition \ref{pro:ProjIsomorphism}) It is immediate from the
definition that $\rho_{i}$ is a poset map, so it suffices to produce
an inverse map. Let $\phi_{i}$ be the map $N\mapsto WN\cap Z$. Since
there is a chain with every edge weakly separated by $i$, $N_{i}\cap W=N_{i}\cap Z$
and $N_{i+1}W=N_{i+1}Z$.

Then for $H$ on $[W,Z]$ we have (by repeated application of the
Dedekind identity)

\begin{eqnarray*}
\phi_{i}\nproj_{i}(H) & = & W(N_{i}H\cap N_{i+1})\cap Z=N_{i}H\cap N_{i+1}W\cap Z\\
 & = & N_{i}H\cap N_{i+1}Z\cap Z=N_{i}H\cap Z=H(N_{i}\cap Z)\\
 & = & H(N_{i}\cap W)=H,\end{eqnarray*}
while for $N$ in $\mathcal{N}_{i}$ we get

\begin{eqnarray*}
\nproj_{i}\phi_{i}(N) & = & N_{i}(WN\cap Z)\cap N_{i+1}=WN\cap ZN_{i}\cap N_{i+1}\\
 & = & WN\cap N_{i+1}\cap\nproj_{i}(Z)=N(N_{i}W\cap N_{i+1})\cap\nproj_{i}(Z)\\
 & = & \rho_{i}(W)N\cap\nproj_{i}(Z),\end{eqnarray*}
and for $N$ between $\nproj_{i}(W)$ and $\nproj_{i}(Z)$ we have
$\nproj_{i}\phi_{i}(N)=N$.
\end{proof}
\begin{note}
Our use of the fact that $N$ is in $\nlat_{i}$ in the proof of Proposition
\ref{pro:ProjIsomorphism} is somewhat subtle: it comes in when we
assume that $WN$ is a subgroup. (Otherwise, $\phi_{i}(N)$ is not
necessarily in $L(G)$.)
\end{note}
\begin{cor}
\label{cor:EdgeProjectsToEdge}If $E\subsetdot F$ is a cover relation
in $L(G)$, then $\rho_{i}(E)\subset\rho_{i}(F)$ is a cover relation
in $\nlat_{i}$.
\end{cor}

\section{Labeling $L(G)$}

Proposition \ref{pro:ProjIsomorphism} and Corollary \ref{cor:EdgeProjectsToEdge}
make it clear how to construct an $EL$-labeling of $L(G)$: label
first by the weak separation labeling, then refine by the modular
labeling in the projection to $\nlat_{i}$. 

More precisely, for each distinct $\nlat_{i}=\nlat_{i}(H)$, let $\lambda^{\nlat_{i}}$
be the modular $EL$-labeling of $\nlat_{i}$. Suppose that $E\subsetdot F$
is an edge in $L(G)$, weakly separated by $N_{i+1}/N_{i}$.

Then label the edge with the pair \[
\lambda(E\subsetdot F)=(i,\quad\lambda^{\nlat_{i}}(\nproj_{i}(E)\subset\nproj_{i}(F)\,).\]
 As is usual, pairs $(i,j)$ are ordered lexicographically.

\begin{thm}
$\lambda$ is an $EL$-labeling of $L(G)$.
\end{thm}
\begin{proof}
Lemma \ref{lem:LMchainGivesPreEL} and Corollary \ref{cor:LMchainGivesPreEL-lex}
tell us that any increasing (lexicographically first) chain on $[W,Z]$
is an extension of the chain $\mathbf{c}=\{C_{0}\subsetdot C_{1}\subsetdot\dots\subsetdot C_{m}\}$,
inductively obtained by taking $C_{0}=W$, and $C_{j+1}=N_{i(j)+1}C_{j}\cap Z$,
where $i(j)$ is the maximal index such that $N_{i(j)}C_{j}\cap Z=C_{j}$. 

Every edge on the interval $[C_{j},C_{j+1}]$ is weakly separated
by $N_{i(j)+1}/N_{i(j)}$, so projects to the same $\nlat_{i}$, and
the modular $EL$-labeling on $\nlat_{i}$ gives a unique increasing
(lexicographically first) chain on $[C_{j},C_{j+1}]$, hence a unique
increasing (lexicographically first) extension of $\mathbf{c}$. 
\end{proof}
\begin{note}
A left modular element in $L$ is also left modular in the dual lattice
$L^{*}$, and Lemma \ref{lem:NlatPreserved} and Proposition \ref{pro:ProjIsomorphism}
say the same thing in $L^{*}$ as in $L$. Thus, we could just as
easily take a chief series $G=N_{0}^{*}\innormal N_{1}^{*}\innormal\dots\innormal N_{k}^{*}=1$,
and label via \[
\lambda_{*}(E\dotsupset F)=(i,\lambda_{*}^{i,\nlat_{i}}(\nproj_{i}^{*}(E)\gtrdot\nproj_{i}^{*}(F)),\]
where $N_{i}^{*}/N_{i+1}^{*}$ weakly separates $E\dotsupset F$ and
$\nproj_{i}^{*}$ is the projection to $[N_{i+1}^{*},N_{i}^{*}]$.
Depending on taste, the resulting $EL$-labeling of the dual lattice
may even seem more natural.
\end{note}

\subsection{Descending chains}

If $E\subsetdot F$ satisfies $E\cap N_{i+1}=1$ and $EN_{i+1}=G$
while $F\cap N_{i}=1$ and $FN_{i}=G$, then $EN_{i+1}=FN_{i+1}=G$
and $E\cap N_{i}=F\cap N_{i}=1$. Thus, $E\subsetdot F$ is separated
by $i$, and thus a chain of complements is a descending chain, labeled
$k-1,\dots,1,0$. By Thévenaz's theorem (Theorem \ref{thm:SubgroupLatticeTopology}),
and since an $EL$-shellable lattice has the homotopy type of a bouquet
of spheres in correspondence to the descending chains, the chains
of complements are exactly the descending chains. (This is also straightforward
to verify by induction.) 

Similarly for the dual labeling $\lambda_{*}$. To summarize:

\begin{prop}
\label{pro:DescChainsAreChainsOfComplements}The descending chains
of both $\lambda$ and $\lambda_{*}$ are exactly the chains of complements
of the chief series used.
\end{prop}
Thus, the labelings we have constructed are the ones conjectured by
Shareshian.

\section*{Acknowledgements}

Thanks to Bruce Sagan, John Shareshian, Hugh Thomas, and the anonymous
referee for their comments and suggestions.

\bibliographystyle{amsplain}
\bibliography{4_Users_paranoia_Documents_Research_Master}

\end{document}